\documentclass[12pt,a4paper,final]{article}

\usepackage{url}
\usepackage[a4paper,portrait]{geometry}
\usepackage{amsmath,amsfonts,latexsym,amssymb,amsthm}

\newcommand{\EM}{\ensuremath}

\makeatletter
\newcommand{\@THMSTYLES}{%
  \newtheoremstyle{bodyrm}
  {3pt}
  {3pt}
  {}
  {}
  {\bfseries\sffamily}
  {.}
  { }
  {}
  \newtheoremstyle{bodyit}
  {3pt}
  {3pt}
  {\itshape}
  {}
  {\bfseries\sffamily}
  {.}
  { }
  {}
}
\newcommand{\THMEN}{%
  \@THMSTYLES
  \theoremstyle{bodyit}
  \newtheorem{thm}{Theorem}[section]%
  \newtheorem{cor}[thm]{Corollary}%
  \newtheorem{prop}[thm]{Proposition}%
  \newtheorem{lem}[thm]{Lemma}%
  \theoremstyle{bodyrm}%
  \newtheorem{defi}[thm]{Definition}%
  \newtheorem{xpl}[thm]{Example}%
  \newtheorem{exo}[thm]{Exercise}%
  \newtheorem{hyp}[thm]{Hypothesis}%
  \newtheorem{eur}[thm]{Heuristics}%
  \newtheorem{pro}[thm]{Problem}%
  \newtheorem{rem}[thm]{Remark}%
  \newtheorem{prp}[thm]{Property}%
}
\newcommand{\THMFR}{%
  \@THMSTYLES
  \theoremstyle{bodyit}
  \newtheorem{thm}{Théorème}[section]%
  \theoremstyle{bodyrm}%
  \newtheorem{defi}[thm]{Définition}%
  \newtheorem{xpl}[thm]{Exemple}%
  \newtheorem{rem}[thm]{Remarque}%
}
\makeatother

\makeatletter
\newcommand{\SMALLSECS}{%
 \renewcommand{\section}{\@startsection%
  {section}
  {1}
  {0em}
  {\baselineskip}
  {0.5\baselineskip}
  {\normalfont\large\bfseries}}
 \renewcommand{\subsection}{\@startsection%
  {subsection}
  {2}
  {0em}
  {\baselineskip}
  {0.25\baselineskip}
  {\normalfont\bfseries}}
}
\makeatother

\makeatletter
\providecommand{\timenow}{\@tempcnta\time
\@tempcntb\@tempcnta
\divide\@tempcntb60
\ifnum10>\@tempcntb0\fi\number\@tempcntb
\multiply\@tempcntb60
\advance\@tempcnta-\@tempcntb
:\ifnum10>\@tempcnta0\fi\number\@tempcnta}
\makeatother

\makeatletter
\newcommand{\versiondetravail}{%
 \renewcommand{\@evenfoot}{%
 \hfil{\tiny\texttt{%
   Version préliminaire, compilée le \today{} à \timenow.}\hfill}}%
 \renewcommand{\@oddfoot}{\@evenfoot}%
}
\makeatother

\newcommand{\dN}{\EM{\mathbb{N}}}

\newcommand{\dP}{\EM{\mathbb{P}}}

\newcommand{\dR}{\EM{\mathbb{R}}}

\newcommand{\rH}{\EM{\mathrm{H}}}

\newcommand{\rL}{\EM{\mathrm{L}}}

\newcommand{\cC}{\EM{\mathcal{C}}}
\newcommand{\cD}{\EM{\mathcal{D}}}

\newcommand{\cF}{\EM{\mathcal{F}}}
\newcommand{\cG}{\EM{\mathcal{G}}}

\newcommand{\cL}{\EM{\mathcal{L}}}

\newcommand{\cN}{\EM{\mathcal{N}}}

\newcommand{\cP}{\EM{\mathcal{P}}}

\newcommand{\cT}{\EM{\mathcal{T}}}

\newcommand{\bF}{\EM{\mathbf{F}}}
\newcommand{\bG}{\EM{\mathbf{G}}}
\newcommand{\bH}{\EM{\mathbf{H}}}

\newcommand{\bK}{\EM{\mathbf{K}}}
\newcommand{\bL}{\EM{\mathbf{L}}}

\newcommand{\bP}{\EM{\mathbf{P}}}

\newcommand{\al}{\alpha}
\newcommand{\be}{\beta}

\newcommand{\de}{\delta}
\newcommand{\ga}{\gamma}

\newcommand{\la}{\lambda}
\newcommand{\La}{\Lambda}

\newcommand{\si}{\sigma}
\newcommand{\Te}{\Theta}

\newcommand{\veps}{\varepsilon}
\newcommand{\vphi}{\varphi}

\newcommand{\p}[4]{{#3}\!\left#1{#4}\right#2} 

\newcommand{\ABS}[1]{\EM{{\left| #1 \right|}}} 
\newcommand{\BRA}[1]{\EM{{\left\{#1\right\}}}} 
\newcommand{\NRM}[1]{\EM{{\left\| #1\right\|}}} 
\newcommand{\PAR}[1]{\EM{{\left(#1\right)}}} 
\newcommand{\SBRA}[1]{\EM{{\left[#1\right]}}} 

\newcommand{\rentf}{\mathbf{Ent}}
\newcommand{\rent}[2]{\p(){\rentf}{#1\,\vert\,#2}}

\newcommand{\Det}[1]{\mathrm{Det}\,}

\newcommand{\GR}{\nabla}

\newcommand{\ADH}[1]{\mathbf{adh}(#1)}

\newcommand{\WH}[1]{\widehat{#1}}

\renewcommand{\leq}{\leqslant}
\renewcommand{\geq}{\geqslant}

\SMALLSECS
\THMEN
\title{On nonparametric maximum likelihood \\
  for a class of stochastic inverse problems}

\author{Djalil~\textsc{Chafa\"\i} \& Jean-Michel~\textsc{Loubes}}

\date{\footnotesize Preprint -- November 2004}

\newcommand{\mykeywords}{%
 Inverse Problems;
 Nonlinear Models;
 Maximum Likelihood;
 EM Algorithm;
 Mixtures of Probability Measures;
 Repeated Measurements Data;
 Longitudinal Data.
}

\newcommand{\mysubjclass}{%
 62G05; 
 34K29. 
}

\begin{document}

\maketitle

\begin{abstract}
  We establish the consistency of a nonparametric maximum likelihood
  estimator for a class of stochastic inverse problems. We proceed by
  embedding the framework into the general settings of early results of
  Pfanzagl related to mixtures \cite{MR89g:62063,Pf2}.
\end{abstract}

{
\footnotesize
 \noindent%
 \textbf{Subject Classification MSC-2000:} \mysubjclass
}

\section*{Introduction}

Let $(S_i,T_i)_{i\in\dN^*}$ be a sequence of i.i.d. random variables with
values in $\dR^p\times\dR_+^n$ and with common law $\mu_S\otimes\mu_T$. Let
$(\veps_i)_{i\in\dN^*}$ be a sequence of i.i.d. standard normal random
variables on $\dR^n$, independent of the preceding sequence. We consider in
the sequel the inverse problem which consists in estimating the law $\mu_S$
given the finite sequence $(Y_i,T_i)_{1\leq i\leq N}$ where
\begin{equation}\label{eq:pb}
  Y_i := f(S_i,T_i)+\si\,\veps_i,
\end{equation}
and where $f:\dR^p\times\dR^n\to\dR^n$ is a known smooth function, which can
be in particular nonlinear in the first variable. The asymptotic is taken in
$N$, and $n$ remains fixed. It is assumed that $\si$ is some known
non-negative variance parameter. We emphasise the fact that in the triplet
$(Y_i,T_i,S_i)$, we observe only the couple $(Y_i,T_i)$, and we are interested
in the estimation of the joint law of the unobserved random variables $S_i$.

In the sequel, $\cL(Z)$ denotes the law of the random variable $Z$. For
example, one has $\cL(S_i,T_i)=\mu_S\otimes\mu_T$. In the same spirit,
$\cL(Z_1\,\vert\,Z_2)$ denotes the conditional law of $Z_1$ given $Z_2$.
Finally, we denote by $\cP(\dR^d)$ the convex set of probability measures on
$\dR^d$ equipped with its Borel $\si$-field and with the $\cC_b(\dR^d,\dR)$
dual topology. We will sometimes denote $S$, $T$, $Y$ for any random variable
with law $\mu_S=\cL(S_1)$, $\mu_T=\cL(T_1)$, and $\mu_Y=\cL(Y_1)$
respectively. Finally, we will denote by $y_i=(y_{i,1},\ldots,y_{i,n})$,
$t_i=(t_{i,1},\ldots,t_{i,n})$ and $s_i=(s_{i,1},\ldots,s_{i,p})$ any
realisation of the random variables $Y_i$, $T_i$ and $S_i$ respectively.

Before starting the mathematical analysis of the problem, let us give briefly
some explanations regarding the notations and the motivations. The random
variables $Y_i=(Y_{i,1},\ldots,Y_{i,n})$ represents the values measured for
individual number $i$ at times $T_i=(T_{i,1},\ldots,T_{i,n})$. The random
variable $S_i$ stands for the individual parameter and the random variable
$\si\,\veps_i$ models the (homoscedastic) random noise which is added to the
possibly nonlinear true value $f(S_i,T_i)$.  This kind of data is known as
\emph{repeated measurements}, or called \emph{longitudinal} since each
individual (from $i=1$ to $i=N$) is observed $n_i:=n$ times and provides a
whole vector of consecutive observations $Y_i=(Y_{i,1},\ldots,Y_{i,n})$
performed at the corresponding individual times
$(T_{i,1},\ldots,T_{i,n})=T_i$. Since $T_i$ is a sequence of measuring times,
one can assume for simplicity that the law $\mu_T$ is a tensor product of
uniform laws on disjoint consecutive compact intervals of the real half line
$\dR_+$. One can think about $\mu_T$ and $n$ as the design of the experiment,
whereas $f$ and $\cL(S_i\,\vert\,T_i)$ and $\cL(\veps_i)$ correspond to the
model chosen for the inverse problem, relating the individual observation
$Y_i$ to the individual parameter $S_i$ and to the individual measuring times
$T_i$. Usually in applications, $f$ is of the form
\begin{equation}\label{eq:split-of-f}
  f(s,T_i) = (q_s(T_{i,1}),\ldots,q_s(T_{i,n})),
\end{equation}
where for any $s$ in $\dR^p$, $q_s:\dR\to\dR$ is a smooth function depending
smoothly on the parameter $s$, for example a linear combination of time
dependent exponentials with coefficients related to $s$. Function $q_s$
represents in such a scheme the true evolution in time of the phenomenon of
interest for an individual of parameter $s$.

Practical applications of models like \eqref{eq:pb} are numerous in signal
transmission, in tomography, in econometrics, in geophysics, etc, cf.
\cite{laref}. Let us give briefly a concrete example in Biology. We consider
the decay of the concentration of a medicine in human blood. One has $p=2$ and
$q_{(A,\al)}(t)=A\exp(-\al t)$ in \eqref{eq:split-of-f}, where $A$ stands for
the quantity of medicine in the blood at time $0$, and where $\al$ stands for
the rate at which the medicine is eliminated. At the beginning of the
experiment, the medicine is given to $N$ independent patients. For patient
number $i$, $n$ measurements $(Y_{i,j})_{1\leq j\leq n}$ of the concentration
of the medicine in blood are made, at times $(t_{i,j})_{1\leq j\leq n}$. One
of the simplest model used in this context is
\begin{equation*}
  Y_{i,j} = q_{(A_{i},\al_i)}(t_{i,j})
  + \si \veps_{i,j},\text{ with }i=1,\ldots,N \text{ and }j=1,\ldots,n.
\end{equation*}
If we state $S_i:=(A_i,\al_i)$, the random variables $S_1,\ldots,S_N$ are
i.i.d. and correspond to the biological specificity of each patient. We are
interested in the estimation of the distribution $\mu_S$ of the common law of
these random variables (population pharmacokinetics). Deconvolution methods
are useless since the required condition $n\to+\infty$ is unrealistic. The
number of observations $n$ for each individual remains small, a few units in
practice. Our framework where the asymptotic is taken on the number of
individuals $N$ is the only mean to perform the estimation of the ``population
law'' $\mu_S$.

A stochastic inverse problem is an inverse problem for which the subject of
the inversion is a probability measure, like in \eqref{eq:pb}. The related
theoretical and applied literature is huge, with many connected components. It
contains in particular deconvolution problems, mixtures models, (non)linear
mixed effects models, (non)linear filtering problems, etc. Even a common
keyword or phrase like our ``stochastic inverse problems'' is most of the time
missing and/or ignored. Therefore, it is quite hard to give a descent state of
the art, but a bit less difficult is to show various natures of a particular
subclass of problems.

We emphasise the fact that \eqref{eq:pb} is not a standard regression problem
since $f$ is \emph{known} whereas the $S_i$ and their law are \emph{unknown}.
Moreover, our problem \eqref{eq:pb} in not of Ibragimov and Hasminskii type
since the $S_i$ are not observed. Notice that when $n$ is very large
deconvolution techniques can give an estimation of each $S_i$. The approach
developed recently in \cite{arnaque} in useless for our problem since we
consider an asymptotic in $N$ and not in $n$.

One of the common difficulties of stochastic inverse problems like
\eqref{eq:pb} lies in the fact that they are ill-posed. The inverse of the
underlying operator is not continuous in general, so that a small perturbation
of the data may induce a large change for the common law of the unobserved
random variable. If the unknown was a function in a Hilbert space instead of a
probability density function, one could try a singular value decomposition
(SVD), following for example Cavalier, Golubev, Picard and Tsybakov in
\cite{Cavalier}.

Several authors have investigated nonparametric maximum likelihood estimation
for stochastic inverse problems, and related Expectation Maximisation (EM, cf.
\cite{MR58:18858}) algorithms. In the context of mixtures, Lindsay showed in
\cite{MR85m:62008a,MR85m:62008b} by using elementary convex analysis that the
fully nonparametric maximum likelihood is achieved by a discrete probability
measure with finite number of atoms related to the sample size, connecting by
this way this kind of problems with convex analysis algorithms (Simplex
algorithm, Fedorov methods, etc). One can find some developments in
\cite{lindsay-monography,MR96h:62075,MR2001a:62001,MR96h:62056}. The
consistency of such estimators was established at least by Pfanzagl in
\cite{MR89g:62063}. In \cite{MR92g:62047}, Schumitzky gave an EM like
algorithm for Lindsay's estimator. In another direction, Eggermont and
Lariccia have developed smoothing techniques for problems involving Fredholm
integral operators, cf. \cite{eggbook} and references therein.

To sum up, our aim in this paper is to estimate $\mu_S$, the common law of the
unobserved i.i.d. random variables $S_i$ in (\ref{eq:pb}), when $\mu_S$
belongs to some class $\cF_S\subset\cP(\dR^p)$. The rest of the paper is
divided as follows. Section \ref{se:NPML} introduces a nonparametric
Likelihood Estimator (NPML) for $\mu_S$, and is devoted to establish its
consistency up to identifiability. Section \ref{se:NPML-algos} presents finite
dimensional and algorithmic approaches to approximate the NPML. Finally, in
Section \ref{se:discus}, various related questions are discussed.

\section{An NPML and its consistency}
\label{se:NPML}

Conditionally on the $S_i$, the $Y_i$ are independent but not identically
distributed, due to the dependence over $T_i$. However, since the individual
observed datum consists in $X_i:=(Y_i,T_i)$, it is quite natural to see $S_i$
as the unique unobserved random variable in the triplet $(Y_i,S_i,T_i)$. The
law $\cL(X_i)=\cL(Y_i,T_i)$ is nothing else but
\begin{equation*}
  \int_{s\in\dR^p}\!\ga_{\si,n}(y-f(s,t))\,d\mu_T(t)\,d\mu_S(s)\,dy,
\end{equation*}
where ``$(y,t)=x$'' and where $\ga_{\si,n}$ is the Gaussian
probability density function on $\dR^n$ given by
$\ga_{\si,n}(u):=(2\pi\si^2)^{-n/2}\,\exp(-\NRM{x}_2^2/2\si^2)$.
Similarly, the law $\cL(Y_i)$ of $Y_i$ is the following mixture
\begin{equation*}
  \SBRA{\int_{s\in\dR^p}\!%
    \int_{t\in\dR_+^n}\!%
    \ga_{\si,n}(y-f(s,t))\,d\mu_T(t)\,d\mu_S(s)}\,dy,
\end{equation*}
where the mixing law is $\mu_S\otimes\mu_T$ and where the mixed family is the
following $f$-deformed Gaussian location family
\begin{equation*}
  \BRA{\ga_{\si,n}(\bullet-f(s,t))\text{ where }(s,t)\in\dR^p\times\dR^n}
  =\ga_{\si,n}*\BRA{\de_{f(s,t)}\text{ where }(s,t)\in\dR^p\times\dR^n}.
\end{equation*}
Assume now that the law $\mu_T$ has a density $\psi$ with respect to the
Lebesgue measure on $\dR_+^n$. Then, one has that the law
$\cL(X_i)=\cL(Y_i,T_i)$ is absolutely continuous with respect to the Lebesgue
measure on $\dR^n\times\dR_+^n$ with probability density function $\bK(\mu_S)$
given by
\begin{equation}\label{eq:def-K}
  \bK(\mu_S)(y,t):=%
  \psi(t)\,\int_{s\in\dR^p}\!\ga_{\si,n}(y-f(s,t))\,d\mu_S(s).
\end{equation}
When $\mu_S$ has density $\vphi$ with respect to Lebesgue's measure on
$\dR^p$, we will denote $\bK(\vphi)$ instead of $\bK(\mu_S)$, viewing by this
way $\bK$ as a linear operator over probability density functions.
\begin{equation*}
  \bK(\vphi)(y,t)=
  \psi(t)\,\int_{s\in\dR^p}\!\ga_{\si,n}(y-f(s,t)) \vphi(s) \,ds. 
\end{equation*}
Here again, the law $\cL(X_i)=\cL(Y_i,T_i)$ is a mixture, with mixing law
$\mu_S$ and mixed family
\begin{equation*}
  \BRA{%
    \psi(t)\,\ga_{\si,n}(\bullet-f(s,t))%
    \text{ with }(s,t)\in\dR^p\times\dR^n}.
\end{equation*}
Notice that $\bK(\mu_S)(y,t)$ is always positive, and thus, $\log\bK(\mu_S)$
always makes sense. The log-likelihood can be expressed by mean of the unknown
law $\mu_S$ as follows
\begin{equation}\label{eq:def-LN}
  \bL_N(\mu_S) := \dP_N \log \bK(\mu_S),
\end{equation}
where $\dP_N$ is the empirical measure of the sample $(X_i)_{1\leq i\leq N}$
defined by
\begin{equation}\label{eq:def-PN}
  \dP_N := \frac{1}{N} \sum_{i=1}^N \de_{(Y_i,T_i)}.
\end{equation}
Notice that we have used above the standard notation $\dP_N F$ to denote the
expectation of function $F$ with respect to probability law $\dP_N$.  When $f$
is of the form \eqref{eq:split-of-f}, the log-likelihood $\bL_N$ defined above
in \eqref{eq:def-LN} reads
\begin{align*}
  \bL_N(\mu_S) %
  &= \frac{1}{N}\,\sum_{i=1}^N %
  \log\int_{s\in\dR^p}\!%
  \exp\PAR{-\frac{1}{2\si^2}\sum_{j=1}^n (Y_{i,j}-q_s(T_{i,j}))^2}\,d\mu_S(s) %
  \quad + C\\
  &=\frac{1}{N}\,\sum_{i=1}^N\sum_{j=1}^n %
  \log\int_{s\in\dR^p} %
  \exp\PAR{-\frac{1}{2\si^2}(Y_{i,j}-q_s(T_{i,j}))^2}\,d\mu_S(s) %
  \quad + C,
\end{align*}
where
\begin{equation*}
  C:=-\frac{n}{2}\log(2\pi\si^2)%
  +\frac{1}{N}\,\sum_{i=1}^N\log\psi(T_{i,1},\ldots,T_{i,n}).
\end{equation*}
The quantity $C$ does not have any effect on the arg-maximum of the
log-likelihood functional $\bL_N$. In particular, the density $\psi$ of
$\mu_T$ does not play a direct role in the NPML \eqref{eq:NPML} below since
one can rewrite $\bL_N$ as follows
\begin{equation*}
  \bL_N(\mu_S) = \dP_N\log\psi + \dP_N\log\bK^\#(\mu_S),
\end{equation*}
where
\begin{equation*}
  (\bK^\#(\mu_S))(y,t)
  := \int_{s\in\dR^p}\!\ga_{\si,n}(y-f(s,t))\,d\mu_S(s).
\end{equation*}
On any set $\cF$, the arg-maximum of $\bL_N$ is equal to the arg-maximum of
$\bL_N^\#$ defined by
\begin{equation*}
  \bL_N^\#(\mu_S) := \dP_N\log\bK^\#(\mu_S).
\end{equation*}
The functional $\bL_N^\#$ does not depend on $\mu_T$ directly, but only
implicitly via the sample $T_1,\ldots,T_N$ throughout $\dP_N$. However, the
law $\mu_T$ plays a role in identifiability, and the good choice of this law
is always a crucial issue.

\begin{defi}[Identifiability] 
  We say that the mixture model \eqref{eq:pb} is \emph{identifiable} if and
  only if $\bK$ is injective, as a map from $\cF_S$ to $\cP(\dR^n)$. Namely,
  for any couple $(\mu,\nu)\in\cF_S\times\cF_S$ with $\nu\neq\mu$, one has
  $\bK(\mu) \neq \bK(\nu)$ in $\cP(\dR^n)$. Similarly, we say that
  $\mu\in\cF_S$ is \emph{identifiable} in \eqref{eq:pb} if and only if
  $\bK(\nu)\neq\bK(\mu)$ in $\cP(\dR^n)$ for any $\nu\in\cF_S$ with
  $\nu\neq\mu$.
\end{defi}

Clearly, the model is identifiable if and only if every element of $\cF_S$ is
identifiable. Identifiability is essential for any estimation issue of the
true mixing law $\mu_S$. This condition is quite difficult to check in great
generality. However, one can find some clues for example in \cite{Concordet}
and references therein. In practice, and when possible, identifiability must
be checked for the particular model considered, and is deeply related to the
properties of function $f$ and to the distribution $\mu_T$ of the observation
times. We are now able to state the following Theorem.

\begin{thm}[Consistency of NPML]\label{th:NPML}
  Assume that $\cF_S\subset\cP(\dR^p)$ is a compact convex subset of a linear
  space, that the model is identifiable, that $\cL(T)=\psi(t)\,dt$, that
  $\mu_S\in\cF_S$, and that for almost all $(y,t)\in\dR^n\times\dR_+^n$, the
  map $\bK(\bullet)(y,t):\cF_S\to\dR$ is continuous. Then, the NPML estimator
  $\WH{\mu_{S,N}}$ given by
  \begin{equation}\label{eq:NPML}
    \WH{\mu_{S,N}} := \arg\max_{\mu\in\cF_S}\bL_N(\mu)
  \end{equation}
  is well defined, unique, and converges almost surely toward $\mu_S$ when $N$
  goes to $+\infty$.
\end{thm}

\begin{proof}
  The random map $\bL_N$ is a.s. continuous from $\cF_S$ to $\dR$ since the
  map $\bK(\bullet)(y,t):\cF_S\to\dR$ is continuous for any
  $(y,t)\in\dR^n\times\dR_+^n$. By linearity and identifiability of $\bK$ and
  strict concavity of the logarithm, the map $\bL_N$ is a.s. strictly concave.
  Thus, it achieves a.s. a unique sup over the compact convex set $\cF_S$. The
  existence and unicity of the estimator $\WH{\mu_{S,N}}$ is therefore proved.
  Finally, thanks to our choice of settings, the desired consistency result
  follows from \cite[Theorem 3.4]{MR89g:62063} and \cite[Section
  5]{MR89g:62063}, since the required hypotheses are fulfilled:
  \begin{itemize}
  \item\textbf{Condition 1}. $\cF_S$ is a compact Hausdorff space, and a
    subset of a linear space.
  \item\textbf{Condition 2}. For almost all %
    $(y_i,t_i)_{1\leq i \leq N}$, the map $\prod_{i=1}^N\bK(\bullet)(y_i,t_i)$
    is continuous on $\cF_S$ for the topology of $\cF_S$.
  \item\textbf{Condition 3}. For almost all %
    $(y,t)\in\dR^n\times\dR_+^n$, the map $\bK(\bullet)(y,t)$ is concave on
    $\cF_S$.
  \end{itemize}
\end{proof}

\begin{rem}\label{rm:NPML-discuss}
  Let us give various remarks about Theorem \ref{th:NPML} and its extensions.
  \begin{enumerate}
  \item\textbf{Identifiability.} Following again \cite{MR89g:62063}, one can
    relax the identifiability of the model to the identifiability of $\mu_S$,
    but it is not really useful in practice since $\mu_S$ is unknown!  For any
    $x:=(y,t)\in\dR^n\times\dR^n$, let us denote by $k_x:\dR^p\to\dR_+$ the
    function $k_x(s):=\ga_{n,\si}(y-f(s,t))$. Let $\cT$ be the biggest open
    subset of $\dR^n$ such that $\psi>0$ over $\cT$. Then, identifiability of
    the model corresponds to a condition on the set of functions
    $\cC:=\BRA{k_x:\dR^p\to\dR_+\text{ with } x\in\dR^n\times\cT}$ appearing
    in the mixture \eqref{eq:def-K}.  Namely, it must separate the elements of
    $\cF_S$. In other words, when $f$ is smooth, the $\cC$ class must be large
    enough to fully characterise any element of $\cF_S$ by duality as a set of
    test functions for a distribution of order zero in the sense of L.
    Schwartz distributions Theory.  Such a necessary and sufficient separation
    condition relies on both $f$ and $\cT$ and can, depending on the
    particular choice of $\cF$, be weaker than the full injectivity of $f$ in
    the first variable when the second runs over $\cT$. Notice that the
    smoothness of $f$ together with its injectivity in the first variable
    induces in general a ``degree of freedom'' requirement on $(n,p)$.  If
    $\cF_S\subset\cD'(K)$ for some compact subset $K$ of $\dR^p$, then $\cC$
    separates the elements of $\mu_S$ as soon as the vector space spanned by
    $\cC$ is dense in $\cC^\infty(K)$ for the uniform topology.
  \item\textbf{Heteroscedasticity.} At least when the elements of $\cF_S$ are
    compactly supported, Theorem \ref{th:NPML} remains true for a class of
    heteroscedastic models of the form
    \begin{equation}\label{eq:hetero}
      Y_i = %
      f(S_i,T_i) %
      + \si\,\veps_i %
      + g(S_i,T_i)\cdot\veps_i,
    \end{equation}
    where $\si>0$ is known, where $g:\dR^p\times\dR^n\to\dR_+^n$ is a smooth
    function and where the dot mark $''\cdot''$ denotes the component-wise
    vectors multiplication. One can also incorporate a matrix between $g$ and
    $\veps_i$. Notice that condition $g\geq0$ ensures that the variance of the
    conditional law is bounded below by $\si^2$ and thus, the mixture makes
    sense. The mixed family is a location-scale $(f,g)$-deformed Gaussian
    family:
    \begin{equation*}
      \BRA{\ga_{(\si^2+g(s,t)^2)^{1/2},\,n}(\bullet-f(s,t))\text{ where }%
        \,(s,t)\in\dR^p\times\dR^n}.
    \end{equation*}
    In concrete applications, it is quite usual to state that $g$ and $f$ are
    co-linear in the heteroscedastic model above, say $g=\si'\,f$, making the
    noise roughly proportional to the measured value.
  \item\textbf{Non Gaussian noise.}  Theorem \ref{th:NPML} remains true when
    the Gaussian law of the noise $\veps_i$ in \eqref{eq:pb} is replaced by an
    absolutely continuous law with respect to the Lebesgue measure on $\dR^n$.
    The related location mixed family is not Gaussian in that case, but this
    does not block the derivation of the consistency of the NPML.
  \item\textbf{Non homogeneity via censure.}  Let $(\mathbf{n}_i)_{i\in\dN^*}$
    be a sequence of i.i.d. random variables independent of
    $(S_i,T_i)_{i\in\dN^*}$, with values in the set $\cN_n$ of subsets of
    $\BRA{1,\ldots,n}$, and with common law
    $p_\kappa:=\dP(\mathbf{n_i}=\kappa)>0$ for any $\kappa\in\cN_n$. Assume
    that for each $i$, one has access only to
    $Z_i:=(Y_{i,j},j\in\mathbf{n}_i)$ instead of the whole vector of
    measurements $Y_i:=(Y_{i,1},\ldots,Y_{i,n})$ itself. Then, the new inverse
    problem corresponds to the new sample
    \begin{equation*}
      ((Z_1,T_1,\mathbf{n}_1),\ldots,(Z_N,T_N,\mathbf{n}_N))
    \end{equation*}
    which is the censored version of the original sample with unobserved
    $S_i$ values
    \begin{equation*}
      ((Y_1,S_i,T_1),\ldots,(Y_N,S_i,T_N)).
    \end{equation*}
    The problem is that the $Z_i$ are not in the same space, but are still
    independent.  Our goal then is to rewrite the problem in a i.i.d
    framework. One method consists in extending the data space to the larger
    direct sum space $E:=\oplus_{\kappa\in\cN_n} E_\kappa$, where $E_\kappa$
    is a copy of $\dR^{\ABS{\kappa}}$ corresponding to the components present
    in $\kappa$, where $\ABS{\kappa}:=\#\kappa$.  It is then easy to write
    down the law of $(Z_i,T_i,\mathbf{n_i})$. Such a model is quite heavy to
    write down but gives rise to a simple extended log-likelihood:
    \begin{equation*}
      \bL_N(\mu_S):= %
      \dP_N\log p_\kappa %
      +\dP_N\log\psi %
      +\dP_N \log %
      \bK_\kappa(\mu_S),
    \end{equation*}
    where for any $\mu\in\cF_S$
    \begin{equation*}
      \bK_\kappa(\mu)(z,t,\kappa) %
      :=\int_{s\in\dR^p}\! %
      \ga_{\si,\ABS{\kappa}}(z-\pi_\kappa(f(s,t)))\,d\mu(s),
    \end{equation*}
    where $\pi_\kappa$ is the projection of $E$ on $E_\kappa$ and where the
    empirical measure $\dP_N$ is now
    \begin{equation*}
      \dP_N:=\frac{1}{N}\sum_{i=1}^N \de_{(Z_i,T_i,\mathbf{n}_i)}.
    \end{equation*}
    The $\dP_N\log p_\kappa +\dP_N\log\psi$ part of the log-likelihood does
    not depend on $\mu_S$, and thus, it does not influence the arg-maximum of
    the log-likelihood and can be safely removed.  For each $i$, the $T_{i,j}$
    involved in the log-likelihood are those with $j\in\mathbf{n_i}$.
    Finally, one can notice that such type of independent censoring does not
    correspond to all realistic censure, since in practice, the $\mathbf{n}_i$
    can depend on the $Y_i$ it self via for example
    \begin{equation*}
      \PAR{\mathrm{I}_{\BRA{Y_{i,1}>\tau}},\ldots,%
          \mathrm{I}_{\BRA{Y_{i,n}>\tau}}}
    \end{equation*}
    where $\tau$ is a detection threshold.    
  \item\textbf{Continuity of the operator.} The continuity assumption on $\bK$
    relies in general on function $f$, on the nature of $\cF_S$, and on the
    law of the noise $\veps_i$, which is Gaussian and homoscedastic here. Some
    concrete examples of $\cF$ are given below.
  \item\textbf{Full extension.} Mixing all the previous extensions is
    delicate. 
  \end{enumerate}
\end{rem}

\begin{xpl}
  Consider for instance the set $\cF_S\subset\cP(\dR^p)$ defined by
  \begin{equation} \label{unifconv}
  \cF_S:=\cF_S^{M,A}:=\BRA{ \vphi(s)\,ds;
    \text{ where }\vphi \in\cC_K^1([0,M]) \text{ and }
    \NRM{\vphi}_{\rL^1}=1,\: \NRM{\ABS{\GR\vphi}}_{\infty}\leq A},
  \end{equation}
  where $K$ is a fixed compact subset of $\dR^p$ and where $M, \: A$ are fixed
  non negative real numbers. Equipped with the $\rL^\infty$ topology, this set
  is a compact convex subset of a linear space, as required by Theorem
  \ref{th:NPML}. Since the underlying mixture model is a ``Gaussian position''
  one, we get for any couple $(\vphi_1,\vphi_2)\in\cF_S\times\cF_S$ and any
  $(y,t)\in\dR^n\times\dR^n$
  \begin{equation*}
    \ABS{\bK(\vphi_1)(y,t)-\bK(\vphi_2)(y,t)} %
    \leq %
    \NRM{\vphi_1-\vphi_2}_\infty \NRM{\psi}_\infty%
    (2\pi\si^2)^{-n/2},
  \end{equation*}
  which gives the $\rL^\infty$ continuity of $\bK(\bullet)(y,t)$ for any
  couple $(y,t)\in\dR^n\times\dR^n$.  Since we deal with a ``Gaussian position
  model'' (homoscedasticity), the operator norm does not depend on $(y,t)$ and
  function $f$ plays not role.  The $\rL^\infty$ a.s.  consistency up to
  identifiability of the NPML follows then from Theorem \ref{th:NPML}.
\end{xpl}

\begin{xpl}
  Consider the set $\cG_S\subset\cP(\dR^p)$ defined by
  \begin{equation*}
    \cG_S:=\cG_S^{A,\al}:=\BRA{\vphi(s)\,ds; %
      \text{ where } \vphi \in \rH^\al(K) %
      \text{ with } %
      \NRM{\vphi}_{\rL^1}=1 %
      \text{ and } \NRM{\vphi}_{H^{\al}}\leq A},
  \end{equation*}
  where $K$ is a fixed compact subset of $\dR^p$, $A$ is a fixed non negative
  real number and $\rH^\al(K)$ is the Sobolev space over the compact $K$.
  Provided that $\al>\frac{1}{2}-\frac{1}{p}$, Rellich-Sobolev embedding
  Theorem yields that $\cF_S$ is a compact convex subset of a linear space for
  the $\rL^2$ topology, cf. \cite{MR56:9247,MR87g:46056}, as required by
  Theorem \ref{th:NPML}.  Since the underlying mixture model is a ``Gaussian
  position'' one, we get for any couple $(\vphi_1,\vphi_2)\in\cF_S\times\cF_S$
  and any $(y,t)\in\dR^n\times\dR^n$
  \begin{equation*}
    \ABS{\bK(\vphi_1)(y,t)-\bK(\vphi_2)(y,t)} %
    \leq %
    \NRM{\vphi_1-\vphi_2}_2 \NRM{\psi}_\infty%
    (4\pi\si^2)^{-n/4},
  \end{equation*}
  which gives the $\rL^2$ continuity of $\bK(\bullet)(y,t)$ for any couple
  $(y,t)\in\dR^n\times\dR^n$. Since we deal with a ``Gaussian position model''
  (homoscedasticity), the operator norm does not depend on $(y,t)$ and
  function $f$ plays not role.  The $\rL^2$ a.s.  consistency up to
  identifiability of the NPML follows then from Theorem \ref{th:NPML}.
\end{xpl}

\section{Algorithms for the NPML}
\label{se:NPML-algos}
\subsection{Finite dimensional approximation}

The first step towards a practical implementation is to transform the maximum
$\WH{\mu_{S,N}}$ of the log-likelihood $\bL_N$ over the whole infinite
dimensional class $\cF_S$ into a maximum $\WH{\mu_{S,N,m}}$ over a finite
dimensional convex subset $\cF_{S,m}$, where $(\cF_{S,m})_{m\in\dN^*}$ is an
exhaustive sequence of subsets of $\cF_S$, i.e. $\ADH{\cup_{m\in\dN^*} \cF_m}
= \cF$.

\begin{thm}\label{th:findim-NPML}
  Assume that $\cF_S$ is a metric space.  Let $(\cF_{S,m})_{m\in\dN^*}$ be an
  exhaustive sequence of finite dimensional closed convex subsets of $\cF_S$.
  Under the assumptions of Theorem \ref{th:NPML}, and for any fixed sample of
  size $N$, the approximated NPML estimator $\WH{\mu_{S,N,m}}$ given by
  \begin{equation}\label{eq:NPML-sieves}
    \WH{\mu_{S,N,m}}%
    := \arg\max_{\mu\in\cF_{S,m}} \bL_N(\mu).
  \end{equation}
  is well defined, unique, and converges toward the NPML $\WH{\mu_{S,N}}$ when
  $m$ goes to $+\infty$.
\end{thm}

\begin{proof}
  We proceed at fixed $N$.  Since $\cF_{S,m}$ is a compact convex subset, the
  approximated NPML estimator $\WH{\mu_{S,N,m}}$ exists, as it was the case
  for the NPML estimator $\WH{\mu_{S,N}}$ in Theorem \ref{th:NPML}. Let us now
  establish the convergence.  By the definition of $\WH{\mu_{S,N,m}}$ and
  $\WH{\mu_{S,N}}$ one has that
  \begin{equation*}
    \bL_N(\WH{\mu_{S,N,m}}) \leq \bL_N(\WH{\mu_{S,N}}).
  \end{equation*}
  In the other hand, there exists a sequence $(\mu_m)_{m\in\dN^*}$ converging
  towards $\WH{\mu_{S,N}}$ in $\cF_S$ and such that $\mu_m\in\cF_{S,m}$ for
  any $m\in\dN^*$. Hence, lower semi continuity of $\bL_N$ induces that, for
  any $\veps>0$, there exists $m_\veps\in\dN^*$ such that for any $m\geq
  m_\veps$,
  \begin{equation*}
    \bL_N(\WH{\mu_{S,N}}) - \veps \leq \bL_N(\mu_m).
  \end{equation*}
  But by definition of $\WH{\mu_{S,N,m}}$ we have
  \begin{equation*}
    \bL_N(\mu_m) \leq \bL_N(\WH{\mu_{S,N,m}}).
  \end{equation*}
  As a result, the following bound holds for any $\veps >0$ and any
  $m>m_{\veps}$
  \begin{equation}\label{eq:sandwiche}
    \bL_N(\WH{\mu_{S,N}}) - \veps \leq \bL_N(\WH{\mu_{S,N,m}})
    \leq \bL_N(\WH{\mu_{S,N}}).
  \end{equation}
  If $\mu^*\in\cF_S$ is an adherence value of the sequence
  $(\WH{\mu_{S,N,m}})_{m\in\dN^*}$, corresponding to the limit point of a
  subsequence $(\WH{\mu_{S,N,m_k}})_{k\in\dN^*}$, then $\mu^*=\WH{\mu_{S,N}}$.
  Namely, if it was not the case, then \eqref{eq:sandwiche} will implies that
  $(\bL(\WH{\mu_{S,N,m_k}}))_{k\in\dN^*}$ converges toward
  $\bL_N(\WH{\mu_{S,N}})$, and thus that $\bL_N(\mu^*)=\bL_N(\WH{\mu_{S,N}})$,
  which contradicts the unicity of $\WH{\mu_{S,N}}$ as a maximum of $\bL_N$
  over $\cF_S$.  Hence, $\WH{\mu_{S,N}}$ is the unique adherence value of the
  sequence $(\WH{\mu_{S,N,m}})_{m\in\dN^*}$, and the compacity of $\cF_S$
  yields finally that $(\WH{\mu_{S,N,m}})_{m\in\dN^*}$ converges towards
  $\WH{\mu_{S,N}}$, which is exactly the desired result.
\end{proof}

\begin{rem}
  The rate of convergence of $(\WH{\mu_{S,N,m}})_{m\in\dN^*}$ towards
  $\WH{\mu_{S,N}}$ when $m$ goes to $+\infty$ depends on the regularity of
  $\cF_{S,m}$ and $\bL_N$.
\end{rem}

\subsection{A Gradient algorithm for log-likelihood maximisation} 
\label{sgrad}

Since for any $m\in\cF_{S,m}$ and any couple $(\mu,\nu)$ in
$\cF_{S,m}\times\cF_{S,m}$,
\begin{equation*}
  \bL_N(\mu)-\bL_N(\nu) = \dP_N \log \frac{ \bK(\mu)}{\bK(\nu)},
\end{equation*}
the sieves log-likelihood estimator $\WH{\mu_{S,N,m}}$ defined in
\eqref{eq:NPML-sieves} can be viewed as the solution of the following
optimisation issue:
\begin{equation}\label{eq:NPML:equiv}
  \text{find } \WH{\mu_{S,N,m}} \text{ such that } %
  \forall \mu \in \cF_{S,m},\ %
  \dP_N \log \frac{ \bK(\mu)}{\bK(\WH{\mu_{S,N,m}})} \leq 0.
\end{equation}
By using the concavity of the objective function, Pfanzagl has proved in
\cite{Pf2} that one may switch, in the definition of the estimator in
\eqref{eq:NPML:equiv}, from the $\log$ function to any other function
$L:\dR_+^*\to\dR$, provided that it is concave, strictly increasing, with
$L(1)=0$.
\begin{equation} \label{eq:Pftrick} 
  \text{find } \WH{\mu_{S,N,m}} \text{ such that } %
  \forall\mu\in\cF_{S,m}, %
  \quad %
  \dP_N L\SBRA{\frac{\bK(\mu)}{\bK(\WH{\mu_{S,N,m}})}} %
  \leq 0.
\end{equation}
As a result defining the estimator for a particular $L$ is enough to get
inequality \eqref{eq:Pftrick} for all ``contrast'' function $L$ satisfying the
previous assumptions. In particular, the estimator $\WH{\mu_{S,N,m}}$ can be
obtained for the special choice $L(t)=t-1$, which corresponds exactly to the
definition of the EM algorithm iteration. Hence, maximising the estimator can
be practically computed via the EM algorithm, while Theorem
\ref{th:findim-NPML} still applies, proving consistency of the estimator. This
invariance in $L$ relies on the ``concavity'' of the model, as explained in
\cite{Pf2}.

\section{Discussion}
\label{se:discus}

\subsection{Heuristics for the NPML in Theorem \ref{th:NPML}}

As usual for maximum log-likelihood, the strong law of large numbers yields
that $(\dP_N)_{N\in\dN^*}$ converges a.s. toward $\bK(\mu_S)$ in $\cP(\dR^n)$.
In other words, $\cL(Y)=\bK(\mu_S)$. Consequently, for any $\mu\in\cF_S$,
$(\bL_N(\mu))_{N\in\dN^*}$ converges toward
\begin{equation*}
  \bL_\infty(\mu):=-\rent{\bK(\mu_S)}{\bK(\mu)}+\bH(\bK(\mu_S)),
\end{equation*}
where $\rent{\bK(\mu_S)}{\bK(\mu)}=\int
(\log\bK(\mu_S)-\log\bK(\mu))\bK(\mu_S)$ is the Kullback-Leibler relative
entropy of $\bK(\mu_S)$ with respect to $\bK(\mu)$ and where
$\bH(\bK(\mu_S))=\bL_\infty(\mu_S)$ is the Shannon entropy of $\bK(\mu_S)$. In
other words, the log-likelihood random functional $\bL_N$ converges toward the
deterministic functional $\bL_\infty$ when $N$ goes to $+\infty$. This
deterministic limit $\bL_\infty$ is the relative entropy functional
$\rent{\bullet}{\bK(\mu_S)}$, up to the additive constant $\bH(\bK(\mu_S))$
which does not play any role for the arg-maximum problem. Since $\bK$ is
injective (identifiability), $\bL_\infty$ is strictly concave with unique
maximum achieved at point $\mu_S$. The NPML estimator replaces the asymptotic
arg-maximum $\mu_S$ with the finite $N$ arg-maximum $\WH{\mu_{S,N}}$. The
non-asymptotic log-likelihood $\bL_N$ is not a relative entropy, but remains
strictly concave. The EM algorithm $\mu_{N,k+1} = \bF_N(\mu_{N,k})$ consists
in approximating $\WH{\mu_{S,N}}$ by finding an entropic lower bound
functional for $\bL_N$ which touches $\bL_N$ at the current step $\mu_{N,k}$.
The EM algorithm in this context can be seen also as a gradient like algorithm
$\mu_{N,k+1}=\mu_{N,k}+\bG_N(\mu_{N,k})$ for the concave functional $\bL_N$,
where $\bG_N$ is the G\^ateau directional derivative of $\bL_N$. It turns out
that this gradient like approach appears as a fixed point iteration
$\mu_{N,k+1}=\bF_N(\mu_{N,k})$ where $\bF_N=\bG_N+\mathrm{Id}$. The fixed
point problem $\bF_N(\mu)=\mu$ corresponds exactly to Bayes rule where the
unknown $\mu_S$ is replaced by the current step $\mu$ and where
$\cL(Y)=\bK(\mu_S)$ is replaced by the first marginal of $\dP_N$. Here again,
$(\bF_N)_{N\in\dN^*}$ converges point-wise toward $\bF_\infty$ which admits
$\mu_S$ as unique fixed point. One of the main feature of EM is the
monotonicity of the objective function $\bL_N$ along the algorithm. The
drawback with such a basic EM approach for nonparametric NPML is the fact that
the support is non increasing along the algorithm.

\subsection{Destruction the log-likelihood concavity for mixtures models}

The log-likelihood of mixtures models is a concave functional of the unknown
mixing probability measure. However, this structure is very sensitive. Lindsay
has showed in \cite{MR85m:62008a} by simply using Minkowski-Caratheodory
Theorem that the fully nonparametric NPML for mixtures models like
\eqref{eq:pb} is achieved by an atomic probability measure with at most $N+1$
atoms. By fully nonparametric, we mean that $\cF_S=\cP(\dR^p)$. This
observation is enough robust to remain valid for heteroscedastic models as in
Remark \ref{rm:NPML-discuss}. Unfortunately, the parametrisation of such
discrete probability measures in terms of weights and support points destroys
the concavity of the log-likelihood objective function $\bL_N$. This lack of
concavity cannot be fixed by the introduction of a stochastic ordering on the
set of discrete probability measures with at most $N+1$ atoms.

\subsection{Semi-parametric estimation}

The convexity structure of the NPML problem is destroyed by the incorporation
of fixed effects estimation. This is typically the case for mixed-effects
models where a linear model structure is imposed to $\mu_S$ and where $\si$ is
unknown in \eqref{eq:pb}. In such cases, the global log-likelihood, seen as a
functional of both random and fixed effects, is not concave and has
potentially many local maxima. The semi-parametric approach developed in
\cite{Pf2} is useless since we do not have a consistent estimator of the fixed
effects regardless of the random effect.

Recall that a typical \emph{mixed effects model} corresponds to some
particular structure (a linear model in general) on the $S_i$ in
\eqref{eq:pb}. Namely, $S_i = \Te{}V_i+\eta_i$, where $V_i$ is an observed
vector of per-individual co-variables (sex, weight, etc), where $\Te$ is an
unknown matrix parameter giving the trend (fixed effect), and where $\eta_i$
is the random effect of unobserved data. In such a model, the
$(V_i)_{i\in\dN^*}$ and the $(\eta_i)_{i\in\dN^*}$ are i.i.d., and the
$\BRA{T_i, V_i, \eta_i, \veps_i,\text{ where } i\in\dN^*}$ are mutually
independent random variables. The goal is then to estimate the $\Te$ matrix
and the common law $\mu_\eta$ of the $(\eta_i)_{i\in\dN^*}$. Such models are
used for example in Biology to let the measurements take into account the
known specificity of each individual while conducting a survey. The pattern,
which is determined by physiological rules is given by the function $f$, while
the specificity of each individual is modelled by the random variables
$(S_i)_{1\leq i \leq N}$. If we write $S_i = \Te{}V_i+m+\eta'_i$ where $m$ is
a fixed parameter to be estimated and where $\eta'_i$ is a centred random
effect, one can first estimate the law of the centred random effect $\eta'$
and then estimate the fixed effects $\Te$ and $m$. However, this approach must
be adapted when the coefficient $\si$ in \eqref{eq:pb} is not known, since it
appears in that case as a new fixed effect to be estimated. We believe that a
semi-parametric extension of our method can be made, providing an estimation
of $(\Te,\mu_\eta)$. The approach presented in \cite{Pf2} does not help since
we do not have a consistent estimator for the fixed effects. Despite the fact
that numerous nonparametric techniques were developed for mixtures models,
the widely used approach in applications of nonlinear mixed effects models is
quite rough and consists in a fully parametric estimation of the first two
moments of the law $\mu_\eta$ of the random effect $\eta$, where it is
arbitrarily assumed that this law is normal or log-normal, cf.
\cite{Mentre94,Mentre97} and \cite{davidian-giltinan} for example. Even if
they speed up the effective computations, such fully parametric approaches are
not satisfactory since the consequences it terms of decision are highly
sensitive to the arbitrarily chosen structure for the random effect law (not
robust).

\subsection{No rates}

To obtain rates of convergence for the maximum likelihood estimator, we
consider a neighbourhood of the true distribution $\mu_S$, defined by the
topology chosen according to fulfils the conditions of Theorem \ref{th:NPML}.
Write $V(\mu_S)$ this neighbourhood, then using compacity there exist a finite
sequence of neighbourhood $V(\mu_k),\: k=1,\dots,r_N$ such that $$
\cF_S-V(\mu_S) \subset \cup_{k=1}^{r_N} V(\mu_k).$$
Hence, finding the rate of
convergence of nonparametric maximum likelihood estimator implies studying the
deviation probability
\begin{align*}
  \bP\PAR{\WH{\mu_{S,N}} \notin V(\mu_S)}
  & \leq \sum_{k=1}^{r_N}  \bP\PAR{\WH{\mu_{S,N}} \in V(\mu_k)} \\
  & \leq \sum_{k=1}^{r_N} \bP\PAR{\sup_{\mu \in V(\mu_k)}\frac{1}{N}
    \sum_{i=1}^N \log
    \left[2\left(1+\frac{(\bK^\#(\mu_S))(X_i)}{(\bK^\#(\mu)(X_i)}
      \right)^{-1}\right] \geq \log \ga}
\end{align*}
for $0<\ga<1$ as it is quoted in \cite{Pf2}. Bounding this deviation
inequality requires two main ingredients. First a bound for the entropy of the
mixture class. Recent works by van der Vaart, see for instance
\cite{vaart2004a} and \cite{vaart2004b}, give upper bounds for the entropy of
such classes and hence provide a control over $r_N$. Second, to conclude,
there is a need for a deviation inequality over the previous empirical
process. Unfortunately, to our concern, concentration bounds in this framework
are very difficult to obtain, preventing further calculations to obtain rates
of convergence. Work in this direction was conducted by van de Geer in
\cite{saramix} but can not be applied in this framework. Thus, it seems rather
difficult to obtain rates of convergence for nonparametric maximum likelihood
estimator using this settings.

\subsection{No sieves}

In order to construct a practical maximum likelihood estimator, one needs to
construct a family of finite dimensional spaces undergoing the assumptions of
Theorem \eqref{th:findim-NPML}. Two main choices are investigated in the
statistical literature, but none fulfils all the needed requirements.

On the one hand, we could consider sieves constructed on log bases. Indeed,
for a basis $(\psi_\la)_{\la\in\La}$ of an Hilbert space, consider for a fixed
integer $m$ the set
\begin{equation*}
 \cF_{S,m}:=\BRA{\vphi \in \cF_S,\text{ s.t. } %
 \log \vphi= \sum_{\la\in\La_m} \be_\la \psi_\la},
\end{equation*}
where $\La_m\subset\La$ with $\ABS{\La_m} \leq m$. If we have taken spline
basis for our initial choice of $\psi_\la$, we get the traditional log-spline
model, well studied by Stone in \cite{Stone90}. Such sets are made of
densities but are not compact for the chosen topology.

On the other hand consider a Multiresolution analysis, see for instance
\cite{Mallat90}, constructed using a wavelet basis, $(\zeta_\la)_{\la\in\La}$.
Hence the finite dimensional sets corresponding to the approximation spaces
are defined by $\cF_{S,m}=\BRA{\vphi=\sum_{\la\in\La_m}\be_\la\zeta_\la}$.
Notice that $\cF_{S,m}$ is a closed convex subset of an Hilbert space.
However, it is not a subset of $\cF_S$, set of the densities. This drawback
appears frequently when estimating densities by wavelet estimators: the
estimate is not a density. This defect, which in standard issues is not
redhibitory, prevents here the use of Theorem \eqref{th:findim-NPML}.

\section*{Conclusion} 

We have shown that the nonparametric maximum likelihood estimator for
\eqref{eq:pb} is consistent. However, the practical construction of usable
sieves in the spirit of Section \ref{se:NPML-algos} is questionable.
Improvements and rates of convergence are difficult to obtain in these
setting. In the case where a large number of observations for each subject are
available, i.e. $n\to+\infty$, the problem can be divided in two sub-issues:
first estimate the random effect and then build a nonparametric estimator of
its density. This point of view is tackled for example in \cite{castillo} or
\cite{Giutys}. However, when there is no hope for more data, in particular
when dealing with medical data for which typically $n$ is less than $5$, we
believe that other types of estimators should be considered.

\providecommand{\bysame}{\leavevmode\hbox to3em{\hrulefill}\thinspace}
\providecommand{\MR}{\relax\ifhmode\unskip\space\fi MR }
\providecommand{\MRhref}[2]{%
  \href{http://www.ams.org/mathscinet-getitem?mr=#1}{#2}
}
\providecommand{\href}[2]{#2}

\begin{center}
  \hrule
\end{center}

{
  \footnotesize
  \noindent
  Djalil~\textsc{Chafa\"\i}.\\
  \textbf{Address:} UMR 181 INRA/ENVT, \'Ecole Nationale V\'et\'erinaire de
  Toulouse, 23 Chemin des Capelles, B.P. 87614, F-31076, Toulouse,
  \textsc{Cedex} 3,
  France.\\
  \textbf{E-mail:} \url{mailto:d.chafai@envt.fr.nospam}\\
  \textbf{Address:} UMR 5583 CNRS/UPS, Institut de Math\'ematiques de
  Toulouse, Universit\'e Paul Sabatier, 118 route de Narbonne, F-31062,
  Toulouse, \textsc{Cedex} 4, France.\\
  \textbf{E-mail:} \url{mailto:chafai@math.ups-tlse.fr.nospam}\\
  \textbf{Web-site:} \url{http://www.lsp.ups-tlse.fr/Chafai/}
  
  \medskip
  
  \noindent
  Jean-Michel~\textsc{Loubes}.\\
  \textbf{Address:} UMR 8628 CNRS/Paris-Sud, B\^atiment 425, D\'epartement de
  Math\'ematiques d'Orsay, Universit\'e d'Orsay Paris XI,
  F-91425, Orsay, \textsc{Cedex}, France.\\
  \textbf{E-mail:} \url{mailto:Jean-Michel.Loubes@math.u-psud.fr.nospam}\\
  \textbf{Web-site:} \url{http://www.math.u-psud.fr/~loubes/}
}

\end{document}